\documentclass[12pt]{article}

\usepackage{amssymb}
\usepackage[all]{xy}

\newtheorem{theorem}{Theorem}[section]
\newtheorem{lemma}[theorem]{Lemma}
\newtheorem{corollary}[theorem]{Corollary}
\newtheorem{definition}[theorem]{Definition}
\newtheorem{proposition}[theorem]{Proposition}

\newtheorem{example}[theorem]{Example}
\newtheorem{remark}[theorem]{Remark}

\newcommand{\Ent}{\mathbb{Z}}
\newcommand{\qed}{\hfill $\Box$ \par \bigskip}

\def\id{{\rm id}}
\def\pf{\noindent {\bf Proof:} }

\def\lu{\rightharpoonup}
\def\sqr#1#2{{\vcenter{\vbox{\hrule height.#2pt\hbox{\vrule width.#2pt
height#1pt \kern#1pt \vrule width.#2pt}\hrule height.#2pt}}}}

\def\m#1{m_{(#1)}}

\def\h#1{h_{(#1)}}

\def\a#1{a_{(#1)}}

\def\l#1{l_{(#1)}}

\def\c#1{c_{(#1)}}

\def\id{{\rm id}}

\def\km{k^{*}}

\def\sym{Sym_{M,n,r}(k)}

\begin{document}

\title{Cocycle twisting of $E(n)$-module algebras and applications to the Brauer group}

\vspace{2cm}

\author{{\normalsize
\begin{tabular}{p{5.7cm}p{1.5cm}p{5cm}}
\hspace{1.1cm} {\large G. Carnovale}\thanks{Corresponding author} &  & \hspace{1.1cm} {\large J. Cuadra} \vspace{4pt}\\
Dipartimento di Matematica &  & Dept.
\'{A}lgebra\\
Pura ed Applicata& & y An\'{a}lisis Matem\'atico \\
via Belzoni 7 & & Universidad de Almer\'{\i}a \\
I-35131 Padua, Italy &  & E-04120 Almer\'{\i}a, Spain \\
email:carnoval@math.unipd.it & & email:jcdiaz@ual.es
\end{tabular}
}}
\date{}
\maketitle

\begin{abstract}
We classify the orbits of coquasi-triangular structures for the
Hopf algebra $E(n)$ under the action of lazy cocycles and the Hopf
automorphism group. This is applied to detect subgroups of the
Brauer group $BQ(k,E(n))$ of $E(n)$ that are isomorphic. For a
triangular structure $R$ on $E(n)$ we prove that the subgroup
$BM(k,E(n),R)$ of $BQ(k,E(n))$ arising from $R$ is isomorphic to
a direct product of $BW(k)$, the Brauer-Wall group of the ground
field $k$, and $Sym_n(k)$, the group of $n \times n$ symmetric
matrices under addition. For a general quasi-triangular structure
$R'$ on $E(n)$ we construct a split short exact sequence having
$BM(k,E(n), R')$ as a middle term and as kernel a central
extension of the group of symmetric matrices of order $r<n$ ($r$
depending on $R'$). We finally describe how the image of the Hopf
automorphism group inside $BQ(k,E(n))$ acts on $Sym_n(k)$.
\end{abstract}



\section*{Introduction}

Since the Brauer group of a field was introduced in 1929, several
versions and generalizations of this important invariant have been
proposed. Apart from its generalization to
commutative rings, which plays an important role in Algebraic
Geometry, one of the most important generalizations was the one
given by Wall in \cite{wall}. He introduced a Brauer group for
$\Ent_2$-graded algebras, nowadays called the Brauer-Wall group,
motivated by the study of quadratic forms and Clifford
algebras, see \cite{wall}, \cite{B}, or \cite{lam}. Wall's
construction was extended to cyclic groups and later to abelian
groups. A great impulse to the theory of graded Brauer groups was
given by Long in 1974 when he constructed a larger Brauer group
for any abelian group, which contains all previous Brauer groups
as subgroups (\cite{LO1}). A step forward was taken in
\cite{Long} by extending his construction to any commutative and
cocommutative Hopf algebra. This new group is today referred to
as the Brauer-Long group. During the late seventies and the
eighties this group was extensively investigated and a
satisfactory knowledge of it was reached. A complete treatment of
the Brauer-Long group is given in the monograph \cite{caenepeel}
by S. Caenepeel. \par \medskip

From a Hopf algebra's point of view the commutativity and
cocommutativity assumption in Long's construction is very
restrictive since many interesting examples of Hopf algebras are
non commutative and/or non cocommutative. This difficulty was
overcome by Caenepeel, Van Oystaeyen and Zhang in \cite{CVZ}
where they constructed a Brauer group for a Hopf algebra with
bijective antipode. The final (for the moment) step was taken in
\cite{VZ} where Pareigis' construction  of the Brauer group of a
symmetric monoidal category \cite{PA} was generalized to the case
of braided monoidal categories. Most of known Brauer groups are
particular cases of this construction and braided monoidal
categories provide the best framework to understand all the
previous constructions. \par \medskip

The Brauer group of a Hopf algebra seems to be much more
complicated than the Brauer-Long group. Evidences of this fact are
provided in \cite{VZ2}, \cite{VZ3}, \cite{VZ4} and \cite{gio} where
subgroups of the Brauer
group of Sweedler's Hopf algebra are investigated, and in
\cite{CG1} where this is also done for the family of Hopf algebras
$H_{\nu}$ introduced by Radford in \cite{acca}. In this paper we
study the Brauer group of the family of Hopf algebras $E(n)$ ($n$
a natural number) introduced in \cite{BDG} and later studied in
\cite{PVO1}, \cite{PVO2}. Each quasi-triangular structure $R$ of
a Hopf algebra $H$ gives rise to a subgroup $BM(k,H,R)$ of the
Brauer group $BQ(k,H)$ of $H$. Dually, each coquasi-triangular
structure $r$ on $H$ induces a subgroup $BC(k, H, r)$ in
$BQ(k,H).$ In \cite{gio}, using cocycle twistings of module
algebras, a strategy to classify all these subgroups was
exhibited. Given a coquasi-triangular Hopf algebra $(H,r)$ one
may consider those cocycles $\sigma:H \otimes H \rightarrow k$
such that the Doi's twisted Hopf algebra
$H^\sigma=\;_{\sigma}\!H_{\sigma^{-1}}$ coincides with $H$ and the
twisted coquasi-triangular structure $r_{\sigma}$ in $H^{\sigma}$ is
easier to manage. Since twisting provides a braided monoidal
equivalence of the categories of comodules and since the Brauer group
is invariant under such an equivalence, $BC(k,H,r) \cong
BC(k,H,r_{\sigma})$. This strategy turned out to be fruitful to
classify all these subgroups for Sweedler's Hopf algebra. \par
\smallskip

In this paper we will use again this technique to classify the
subgroups corresponding to the (co)quasi-triangular structures of
$E(n)$ (Theorem \ref{classification} and Corollary \ref{classi}).
We will focus on those cocycles $\sigma$ for which
$H^{\sigma}=H$. These are called lazy cocycles and they form a
group, denoted by $Z_L(H)$, which acts on the set of
coquasi-triangular structures of $H$ (\cite{lazy}). For two coquasi-triangular
structures on $H$, one is the twist of the other if and only if
they are in the same $Z_L(H)$-orbit. In Section 2 we describe the
orbits of coquasi-triangular structures on $E(n)$ under the
action of $Z_L(E(n))$. Coquasi-triangular structures on $E(n)$ are
parametrized by matrices in $M_n(k)$. We prove in Proposition
\ref{orbits} that for $A \in M_n(k)$ the corresponding
coquasi-triangular structure $r_A$ is a twist of $r_B$ for some
$B \in M_n(k)$ if and only if $A-B$ is a symmetric matrix. As a
consequence all $Z_L(E(n))$-orbits are parametrized by
skew-symmetric matrices (Corollary \ref{skew}) and cotriangular
structures lie in the same orbit as $r_0$ where $r_0$ is the
nontrivial normalized bicharacter of $k\Ent_2$ (Corollary
\ref{symm}). The abelian group $BC(k,E(n),r_0)$ is computed in
the first part of Section 3. Since $E(n)$ is self-dual, the
computation of this group is carried out from a dual point of
view, by computing $BM(k,E(n),R_0)$ where $R_0$ is the nontrivial
triangular structure in $k\Ent_2$. The injection of $k\Ent_2$ into
$E(n)$ is a quasi-triangular map and it induces a split
group homomorphism from $BM(k,E(n),R_0)$ to $BW(k)$, the
Brauer-Wall group of $k$. The kernel of this homomorphism is
shown to be isomorphic $Sym_n(k)$, the group of symmetric matrices of order
$n$ over $k$ under addition. Hence all subgroups $BC(k,E(n),r)$
and $BM(k,E(n),R)$ are described for $r$ and $R$ cotriangular and
triangular respectively. The above kernel is determined by
studying the action of $E(n)$ on central simple algebras and its
interaction with Clifford-type algebras, which is done in Section
3. The non-triangular case is more difficult to manage but we
are able to construct in Theorem \ref{reduction} a split short
exact sequence whose middle term is $BM(k, E(n),R_A)$ and whose
left hand term is described as a nonabelian central extension of
$(Sym_r(k),+)$ with $r<\left[\frac{n}{2}\right]$ by
$(M_{n-r,r}(k),+)$, the group of $(n-r)\times r$ matrices with
coefficients in $k$ under addition.
\par \medskip

Finally we compute a subgroup of $BQ(k, E(n))$ that arises from a
quotient of the Hopf automorphism group of $E(n)$. There is a
homomorphism from the Hopf automorphism group of a Hopf algebra
into its Brauer group. Through this homomorphism $GL_n(k)/\Ent_2$
is embedded into $BQ(k, E(n))$. The action by conjugation of the
former group stabilizes $BM(k, E(n), r_0)$ and we construct a
copy of the semidirect product $(GL_n(k)/\Ent_2)\ltimes Sym_n(k)$ inside
$BQ(k, E(n))$.

\section{Notation and preliminaries}\label{one}

Throughout this paper $k$ will denote a fixed ground field of
characteristic different from $2$ and its group of invertible
elements will be denoted by $k^*$. Unless otherwise stated, all
vector spaces, algebras, tensor products, etc will be over $k$.
For vector spaces $V$ and $W$, $\tau:V \otimes W \rightarrow W
\otimes V$ stands for the usual flip map. By $H$ we will denote a
finite dimensional Hopf algebra with bijective antipode $S$. For
general facts on the theory of Hopf algebras the reader is
referred to \cite{LR} and \cite{Mon}. \par
\medskip

\subsection{The Brauer group of a Hopf algebra}

In this paragraph we briefly recall from \cite{CVZ} and
\cite{CVZ2} the construction of the Brauer group of a Hopf
algebra and record some of its properties needed in the sequel.
Let $R=\sum R^{(1)} \otimes R^{(2)} \in H \otimes H$ be a
quasi-triangular structure on $H$. The {\it braided product} $A \#
B$ of two $H$-module algebras $A,B$ is the $H$-module algebra
defined as follows: as an $H$-module $A \# B=A \otimes B$, with
multiplication given by
$$(a \# b)(a' \# b')=\sum a(R^{(2)}\cdot a')\# (R^{(1)} \cdot b)b',$$
for all $a,a' \in A,b,b' \in B.$ The {\it $H$-opposite algebra} of
$A$, denoted by $\overline{A},$ is equal to $A$ as a $H$-module
but endowed with  multiplication: $aa'=\sum (R^{(2)}\cdot
a')(R^{(1)} \cdot a)$ for all $a,a' \in A$. The endomorphism
algebra $End(M)$ of a left $H$-module $M$ is an $H$-module algebra
with the $H$-module structure defined by
$$(h \cdot f)(m)=\sum h_{(1)} \cdot f(S(h_{(2)}) \cdot m).$$
Similarly, $End(M)^{op}$ equipped with the $H$-module structure:
$$(h \cdot f)(m)=\sum h_{(2)} \cdot f(S^{-1}(h_{(1)}) \cdot m)$$
becomes a left $H$-module algebra. A finite dimensional $H$-module
algebra $A$ is called {\it $H$-Azumaya} if the following
$H$-module algebra maps are isomorphisms:
$$\begin{array}{ll}
F: A \# \overline{A} \rightarrow End(A),\ F(a \# \bar{b})(c)=\sum
a(R^{(2)} \cdot c)(R^{(1)} \cdot b), \smallskip \\
G: \overline{A} \# A \rightarrow End(A)^{op},\ G(\bar{a} \# b
)(c)=\sum (R^{(2)} \cdot a)(R^{(1)} \cdot c)b.
\end{array}$$

Two Azumaya $H$-module algebras $A,B$ are called {\it Brauer
equivalent}, denoted by $A \sim B$, if there are finite
dimensional $H$-modules $M,N$ such that $A \# End(M) \cong B \#
End(N)$ as $H$-module algebras. The relation $\sim$ is an
equivalence relation on the set $Az(H)$ of isomorphism classes of
$H$-Azumaya module algebras and the quotient set
$BM(k,H,R)=Az(H)/\sim$ is a group under the product $[A][B]=[A \#
B]$, with identity element $[End(M)]$ for $M$ a finite
dimensional $H$-module, and inverse $[A]^{-1}=\left[\overline{A}\right]$.
Starting with a coquasi-triangular structure $r$ on $H$, a dual
construction holds. The corresponding Brauer group is denoted by
$BC(k,H,r)$. The duality yields an isomorphism $BM(k,H, R) \cong
BC(k,H^{*},R^*)$ where $R^*$ is the coquasi-triangular structure
on $H^*$ induced by $R$. The Brauer group $BM(k,D(H),{\mathcal
R})$ of the Drinfeld double of $H$ with its canonical
quasi-triangular structure ${\mathcal R}$ is denoted by
$BQ(k,H).$ Identifying the monoidal categories $_H{\mathcal
YD}^H$ of Yetter-Drinfeld $H$-modules and $_{D(H)}{\mathcal M}$,
a $D(H)$-module algebra is a Yetter-Drinfeld $H$-module algebra.
\par \smallskip

If $(H,R)$ is a quasi-triangular Hopf algebra, then $BM(k,H,R)$ is
a subgroup of $BQ(k,H)$. It consists of those classes $[A] \in
BQ(k,H)$ having a representative $A$ for which the right
$H^{op}$-comodule algebra structure is of the form: $$\rho:A
\rightarrow A \otimes H^{op},\ a \mapsto \sum (R^{(2)} \cdot a)
\otimes R^{(1)}.$$

Similarly, if $(H,r)$ is a coquasi-triangular Hopf algebra, then
$BC(k,H,r)$ is a subgroup of $BQ(k,H)$. A class $[A] \in BQ(k,H)$
belongs to $BC(k,H,r)$ if it has a representative $A$ such that
the left $H$-module structure is of the form: $$\cdot: H \otimes A
\rightarrow A,\ h \cdot a = \sum r(h \otimes a_{(1)})a_{(0)}.$$

\subsection{The Hopf algebras $E(n)$}

Let $n$ be a natural number and let $E(n)$ be the Hopf algebra
over $k$ generated by $c$ and $x_i$ for $i=1,\,\ldots,\,n$ with
relations
$$c^2=1,\quad x_i^2=0, \quad cx_i=-x_ic,\quad x_ix_j=-x_jx_i, \quad for
\ 1\leq i<j \leq n,$$
coproduct
$$\Delta(c)=c\otimes c,\quad \Delta(x_i)=1\otimes x_i+x_i\otimes c,$$
and antipode
$$S(c)=c,\quad S(x_i)=cx_i.$$
Each $E(n)$ is a quasi-triangular Hopf algebra. Its
quasi-triangular structures were described by Panaite and Van
Oystaeyen in \cite{PVO1}. Notation and terminology throughout
this paragraph is that of this reference. Since $E(n)\cong
E(n)^*$, we have that $E(n)$ is also coquasi-triangular. An
isomorphism $\phi\colon E(n)\to E(n)^*$ is defined by
$\phi(1)=1^*+c^*$, $\phi(x_j)=x_j^*+(cx_j)^*$ and
$\phi(c)=1^*-c^*$. The quasi-triangular structures are
parametrized by matrices in $M_n(k)$ and they are given as
follows: for a matrix $A=(a_{ij}) \in M_n(k)$ and for $s$-tuples
$P,F$ of increasing elements in $\{1,\ldots,n\}$ we define
$|P|=|F|=s$ and $x_P$ as the product of the $x_j$'s whose index
belongs to $P$, taken in increasing order. Any bijective map
$\eta:F \rightarrow F$ may be identified with an element of the
symmetric group $S_s$. Let $sign(\eta)$ denote the signature of
$\eta$. If $P=\emptyset$ then we take $F=\emptyset$ and
$sig(\eta)=1$. Finally, by $a_{P,\eta(F)}$ we denote the product
$a_{p_1,f_{\eta(1)}}\cdots\,a_{p_s,f_{\eta(s)}}$. For
$P=\emptyset$ we define $a_{P,\eta(F)}:=1$. The $R$-matrix
corresponding to $A$ is
$$\begin{array}{ll}
R_A= &
{1\over2}\sum_P(-1)^{{|P|(|P|-1)}\over2}\sum_{F,\;|F|=|P|,\;
\eta\in S_{|P|}}sign(\eta)a_{P,\eta(F)}\bigl(x_P\otimes
c^{|P|}x_F \vspace{4pt} \\
 & +cx_{P}\otimes c^{|P|}x_F+x_P\otimes c^{|P|+1}x_F-cx_{P}\otimes
c^{|P|+1}x_F\bigr).
\end{array}$$
In particular, $\frac{1}{2}a_{ij}$ is the coefficient of
$(x_i\otimes cx_j+ cx_i\otimes cx_j+x_i\otimes x_j-cx_i\otimes
x_j)$. The following proposition strengthens \cite[Proposition
7]{PVO1}.

\begin{proposition}
The quasi-triangular structure $R_A$ is triangular if and only if
$A$ is symmetric.
\end{proposition}
\pf If $(H,R)$ is a quasi-triangular Hopf algebra, then
$(H^{op},R^{-1})$ and $(H^{op},\tau R)$ are quasi-triangular Hopf
algebras. Since $E(n)\cong E(n)^{op}$, we have that
$(E(n),R_A^{-1})$ and $(E(n),\tau R_A)$ are quasi-triangular Hopf
algebras. The quasi-triangular structure of $E(n)^{op}$
corresponding to $A^t$, the transpose matrix of $A$, is $\tau
R_A$. On the other hand, it is well-known that
$R^{-1}=(S\otimes\id)(R)$ and a simple computation shows that
$(S\otimes \id)(R_A)$ is the quasi-triangular structure of
$E(n)^{op}$ corresponding to $A$. Hence, $\tau R_A=R_A^{-1}$ if
and only if $A$ is a symmetric matrix. \qed

By duality, the coquasi-triangular structures of $E(n)$ are
parametrized by matrices in $M_{n}(k)$ and they can be obtained 
applying the isomorphism $\phi\otimes\phi$ to the quasi-triangular
structures. By direct computation one obtains:
$$\begin{array}{ll}
r_A=& \sum_P (-1)^{{|P|(|P|-1)}\over2}\sum_{F, |F|=|P|, \eta\in
S_{|P|}}sign(\eta)a_{P,\eta(F)}\bigl((x_P)^*\otimes (x_F)^* \vspace{4pt} \\
 & +(cx_{P})^*\otimes (x_F)^*+(-1)^{|P|}(x_P)^*\otimes (cx_F)^*-(-1)^{|P|}
 (cx_{P})^*\otimes (cx_F)^*\bigr).
\end{array}$$
In particular, $r_A(x_i\otimes x_j)=a_{ij}.$ Hence $A$ is nothing
but the matrix of the bilinear form determined by the restriction
of $r_A$ to the span of the $x_j$'s. In other words, the matrix
whose entries are the values of a coquasi-triangular structure $r$
on $x_i\otimes x_j$ for $i,j=1,\,\ldots,\,n$ uniquely determines
$r$. As a consequence of the above proposition, $r_A$ is
cotriangular if and only if $A$ is symmetric. Observe that this
does not mean that the form $r_A$ is symmetric on $E(n)\otimes
E(n)$. For instance, as $$r(x_i\otimes x_j)=r(cx_i\otimes
x_j)=-r(x_i\otimes cx_j)= r(cx_i\otimes cx_j)$$ the form $r_A$
cannot be symmetric if $a_{ii}\not=0$ for some $i$. \par
\smallskip

Finally we recall from \cite[Lemma 1]{PVO2} that the group
$Aut_{Hopf}(E(n))$ of Hopf automorphisms of $E(n)$ is isomorphic
to $GL_n(k)$. Any element $f \in Aut_{Hopf}(E(n))$ is of the form
$$f(c)=c, \qquad f(x_i)=\sum_{j=1}^n t_{ij}x_j \quad \forall
i=1,\ldots ,n,$$ for $T_{f}=(t_{ij}) \in GL_n(k)$.

\section{Cocycle twisting}\label{twistings}

Recall that a {\it left $2$-cocycle} is a convolution invertible
map $\sigma:H \otimes H \rightarrow k$ satisfying:
$$\sum \sigma(g_{(1)} \otimes h_{(1)})\sigma(g_{(2)}h_{(2)}
\otimes m)= \sum \sigma(h_{(1)} \otimes m_{(1)})\sigma(g \otimes
h_{(2)}m_{(2)}),$$ for all $g,h,m \in H$. When $\sigma(h \otimes
1)=\sigma(1 \otimes h)=\varepsilon(h)1,$ for all $h \in H$, the
cocycle is said to be {\it normalized}. A new Hopf algebra
$H^{\sigma}$, called the {\it $\sigma$-twist} of $H$, can be
associated to $\sigma$ (\cite{doitake}). As a coalgebra $H^{\sigma}=H$ but with
multiplication defined by $$h \cdot_{\sigma} h'=\sum
\sigma(h_{(1)} \otimes
h'_{(1)})h_{(2)}h'_{(2)}\sigma^{-1}(h_{(3)} \otimes h'_{(3)})$$
for all $h,h'\in H$. If $(H,r)$ is coquasi-triangular, then
$(H^{\sigma},r_{\sigma})$ is coquasi-triangular where
$r_{\sigma}=\sigma \tau * r * \sigma^{-1}$. It is known that
${\mathcal M}^H$ and ${\mathcal M}^{H^{\sigma}}$ are equivalent as
braided monoidal categories. As a consequence, the Brauer groups
$BC(k,H,r)$ and $BC(k,H^{\sigma},r_{\sigma})$ are isomorphic. In
this section we study those $2$-cocycles $\sigma$ of $E(n)$ such
that $E(n)^{\sigma}=E(n)$ and we describe how the
coquasi-triangular structures of $E(n)$ transform under this
twist procedure. This will allow us to detect isomorphism
classes of groups of type $BC(k,E(n),r_A)$ for a
coquasi-triangular structure $r_A$ on $E(n)$. \par \smallskip

\begin{definition} A left $2$-cocycle $\sigma$ is called
{\em lazy} if for all $h,l \in H$
\begin{equation}\label{condition}
\sum\sigma(\h1\otimes\l1)\h2\l2=\sum \h1\l1\sigma(\h2\otimes\l2).
\end{equation}
\end{definition}

It is shown in \cite{lazy} that lazy cocycles are also right
cocycles and that they form a group $Z_L(H)$ under the
convolution product. This group acts on the set ${\mathcal U}$ of
coquasi-triangular structures of $H$ by $\sigma\cdot
r=(\sigma\tau)*r*\sigma^{-1}$. It is well-known that the group
$Aut_{Hopf}(H)$ of Hopf automorphisms acts on cocycles and on
coquasi-triangular structures by
$\alpha\cdot\sigma=\sigma\circ(\alpha^{-1}\otimes\alpha^{-1})$ and
$\alpha\cdot r=r\circ(\alpha^{-1}\otimes\alpha^{-1})$. The
semidirect product ${\mathcal H}:=Z_L(H)\rtimes Aut_{Hopf}(H)$
acts again on $\cal U$. \par \smallskip

We will analyze these orbits when $H=E(n)$. We first describe
$Z_L(E(n))$. To detect orbits, it is enough to describe normalized
lazy cocycles. Indeed, if $\sigma$ is a lazy cocycle which is not
normalized, then $\sigma^{-1}(1\otimes 1)\sigma$ is a normalized
lazy cocycle (see \cite[Page 231]{caenepeel} for the normalization
condition) and it acts as $\sigma$ on $\cal U$. Therefore we
shall focus on normalized lazy cocycles. \par

\begin{lemma}\label{necessary}Let $\sigma$ be a normalized lazy cocycle for $E(n)$. Then 
$\sigma(c\otimes c)=1$ and $\sigma(c\otimes
x_j)=\sigma(x_j\otimes c)=0$.
\end{lemma}
\pf Condition (\ref{condition}) for $h=c$, $l=x_j$ forces
$\sigma(c\otimes c)=1$ and $\sigma(c\otimes x_j)=0$. Similarly,
for $h=x_j$, $l=c$ condition (\ref{condition}) forces
$\sigma(x_j\otimes c)=0.$\qed
\begin{lemma}\label{sufficient}For every choice of scalars $m_{ij}$
 for $1\le i\le j\le n$ there exists a normalized
 lazy cocycle for $E(n)$ with 
$\omega(x_i\otimes x_j)=m_{ij}$  if $i\le j$ and zero if $i>j$. 
\end{lemma}
\pf By the classification of $E(n)$-cleft extensions of $k$ in
\cite{PVO2}, cohomology classes of normalized cocycles are parametrized by
$(\alpha,\gamma, M)$ where 
$\alpha\in\km$, $\gamma\in k^n$ and $M$ is an 
upper triangular $n\times n$ matrix with entries in $k$. They are
given by:
$$\omega(c\otimes c)=\alpha,\quad \omega(c\otimes x_i)=0,\quad
\omega(x_i\otimes c)=\gamma_i$$ for every $i=1,\,\ldots,\,n$,  
$$\omega(x_i\otimes x_j)=\cases{m_{ij}& if $i\le j$;\cr
0& if $i>j$;\cr}\quad 
$$
$$\omega(c\otimes x_P)=\omega(x_P\otimes c)=\omega(c\otimes
cx_P)=\omega(cx_P\otimes c)=0$$ for every $P$ with $|P|>1$;
$$\omega(x_P\otimes x_Q)=0$$ whenever $|P|\neq|Q|$, together with
recurrence relations.  Let us point out that  we slightly changed
notation in \cite{PVO2} for later convenience.\par
\smallskip
If $h=1$ or if $k=1$, or if $h=k=c$ condition (\ref{condition}) is
always satisfied. By Lemma \ref{necessary} we know that any lazy
cocycle in this family should correspond to a triple with 
$\alpha=1$ and $\gamma=0$. Let $\omega=\omega(M)$ be the cocycle
corresponding to $\alpha=1$, $\gamma=0$ and $M$. For such a cocycle 
condition (\ref{condition})
holds also for $h=x_i$ and $k=x_j$. Besides, the recurrence relations
in \cite{PVO2} become
$$\omega(cx_P\otimes x_Q)=\omega(x_P\otimes
x_Q)=(-1)^{|P|}\omega(x_P\otimes cx_Q)=(-1)^{|P|}\omega(cx_P\otimes
cx_Q)$$ and, for $|Q|=|P|+1$
$$\omega(x_ix_P\otimes x_Q)=\sum_{j=1}^{|Q|}(-1)^{|Q|-j}
m_{iq_j}\omega(x_P\otimes
x_{Q\setminus\{q_j\}}).$$
The elements of $Z_L(E(n))$ are precisely those cocycles which do not
change the product in $E(n)$ when we apply Doi's twisting
procedure, i.e.,
$$h\cdot_{\omega} l=\sum \omega(\h1\otimes\l1)\h2\l2\omega^{-1}(\h3\otimes\l3)=hl$$
for all $h,l \in E(n)$.  
It is not hard to verify that $c\cdot_{\omega}(c^ax_P)=cc^ax_P$ and $x_j\cdot_{\omega}(c^ax_P)=x_jc^ax_P$ for
every $j=1,\ldots, n$, every $n$-tuple $P$ and every $a=0,1$. The assertion
follows by induction and associativity of the product
$\cdot_{\omega}$.\qed


\begin{proposition}\label{orbits}
Two coquasi-triangular structures $r_A$ and $r_B$ on $E(n)$ are in
the same $Z_L(E(n))$-orbit if and only if $A-B$ is a symmetric
matrix.
\end{proposition}
\pf Let $r_A$ be the coquasi-triangular structure corresponding to
the matrix $A$ and let $\sigma$ be a normalized lazy cocycle. Then
$\sigma \cdot r_A=\sigma\tau*r_A*\sigma^{-1}=r_B$ for some matrix
$B$, which is completely determined by:
$b_{ij}=\sigma\tau*r_A*\sigma^{-1}(x_i\otimes x_j)$. By direct
computation we obtain:
$$b_{ij}=-\sigma(x_j\otimes x_i)+a_{ij}+\sigma^{-1}(x_i\otimes x_j)
=a_{ij}-\sigma(x_j\otimes x_i)-\sigma(x_i\otimes x_j).$$ Then
$B=A-L$ where $L$ is the {\em symmetric} matrix with coefficients
$l_{ij}=\sigma(x_i\otimes x_j)+\sigma(x_j\otimes x_i)$. Hence, if
$r_B\in Z_L(H)\cdot r_A$, then $B-A\in Sym_n(k)$. Conversely, let
$A,B \in M_{n}(k)$ be such that $L=A-B\in Sym_n(k)$. By Lemma
\ref{sufficient}  there exists a lazy
cocycle $\omega$ such that $\omega(x_i\otimes
x_j)+\omega(x_j\otimes x_i)=l_{ij}$.\qed
%
%

\begin{remark} {\normalfont The classification of cleft extensions
of $E(n)$ in \cite{PVO2} yields a classification of $2$-cocycles
up to cohomologous cocycles. For a $2$-cocycle $\sigma$, a
cocycle is cohomologous to $\sigma$ if 
it is of the form $$\sigma^{\theta}(h\otimes
l)=\sum
\theta(\h1)\theta(\m2)\sigma(\h2\otimes\m2)\theta^{-1}(\h3\m3)$$
for every $h,\,m\in H$ and
for a convolution invertible
map $\theta\colon H\to k$.
However, if $\sigma$ is lazy,
$\sigma^{\theta}$ need not be lazy and viceversa. Therefore, in order to
describe $Z_L(E(n))$ completely we would need to compute {\em all}
$2$-cocycles first.}
\end{remark}

\begin{corollary}\label{symm}
All cotriangular structures on $E(n)$ form a unique
$Z_L(E(n))$-orbit in $\cal U$.
\end{corollary}
\pf We have seen that a coquasi-triangular structure is
cotriangular if and only if its corresponding matrix is
symmetric, hence if and only if it lies in the orbit of the trivial
matrix.\qed

\begin{remark}{\normalfont Since the action of a lazy cocycle
preserves cotriangularity and since $r_0$ is cotriangular we get
a new proof of the fact that if $A$ is symmetric $r_A$ is
cotriangular.}
\end{remark}

\begin{corollary}\label{skew}
The $Z_L(E(n))$-orbits on $\cal U$ are parametrized by
skew-symmetric matrices.
\end{corollary}
\pf They are parametrized by the skew-symmetric part of the
corresponding matrix. \qed

\begin{theorem}\label{classification}
The ${\mathcal H}$-orbits on $\cal U$ are parametrized by $n\times
n$ matrices of the form
\begin{equation}\label{matrici}
J_l:=\pmatrix{0&I_l&0\cr -I_l&0&0\cr 0&0&0}
\end{equation}
where $I_l$ is the $l\times l$ identity matrix for some $l$ with
$0\le l\le\left[{n\over2}\right]$.
\end{theorem}
\pf We know that $Aut_{Hopf}(E(n))\cong GL_n(k)$. By \cite[Proposition
3]{PVO2} it acts on
$R_A$ (and, dually, on $r_A$) by $\alpha_M\cdot R_A=R_{M^tAM}$,
where $\alpha_M$ is the automorphism associated to the matrix
$M\in GL_n(k)$. Hence the ${\mathcal H}$-orbits are
parametrized by congruence classes of skew-symmetric matrices. By
\cite[Theorem XV.8.1]{Lang} any alternating form is equivalent to
(one and only one) form associated to a matrix of type
$($\ref{matrici}$)$. \qed

\begin{corollary}\label{classi}
The possible non isomorphic Brauer groups $BC(k,\,E(n),\,r_A)$
for $A \in M_n(k)$ are at most $\left[{n\over2}\right]+1$. Every
group $BC(k,E(n), r_A)$ is isomorphic to a group of type $BC(k,
E(n), r_{J_l})$. More precisely, $BC(k,E(n), r_A)$ is isomorphic
to that $BC(k, E(n),R_{J_l})$ for which there exist $L\in
GL_n(k)$ and $C\in M_n(k)$ such that $C=\;^tLJ_lL$ and $C-A$ is
symmetric. In particular, the Brauer groups $BC(k,\,E(n),\,r_A)$
corresponding to cotriangular structures are all isomorphic to
$BC(k,E(n),r_0)$.
\end{corollary}

\begin{remark}\label{dual}{\normalfont
Since $E(n)$ is self-dual, $BM(k,E(n),R_A)$ is isomorphic to
\newline $BC(k,E(n),r_A)$. A result similar to Corollary \ref{classi}
 holds for the groups
$BM(k,E(n),R_A).$}
\end{remark}

\begin{example}
{\normalfont i) If $n=1$ then $E(n)$ is equal to Sweedler's Hopf
algebra $H_4$. There is only one $\cal H$-orbit because
$span(x_j)$ is one-dimensional. All matrices are symmetric so all
forms are cotriangular.  This case was studied in \cite{gio} and
its dual version appears in \cite[Example 7.1.4.2]{gelaki}. The
Brauer group $BM(k,\,H_4,\,R_0)\cong BC(k,\,H_4,\,r_0)$ was
computed in \cite{VZ3}. \par \smallskip

ii) For $n=2$ there are two $\mathcal{H}$-orbits: the symmetric
matrices, corresponding to cotriangular structures, and the
non-symmetric matrices, with representative $$\pmatrix{0&1\cr
-1&0}.$$

iii) If $n=3$ again $\mathcal{H}$ acts on $\mathcal{U}$ with two
orbits, corresponding to the symmetric matrices and the
non-symmetric matrices. The second orbit is represented by the
matrix} $$\pmatrix{0&1&0\cr -1&0&0\cr 0&0&0}.$$
\end{example}
\par \medskip

We end this section by introducing some special lazy cocycles and
their corresponding $E(n)$-cleft extensions of $k$. Let $\omega$ be
the lazy cocycle corresponding to the upper triangular matrix $M$ as
in Lemma \ref{sufficient} and  
let  $\theta=\varepsilon-\sum_{1\le i<j\le n}m_{ij}(x_ix_j)^*\in
 E(n)^*$.
We consider the cocycle 
$\sigma:=\omega^\theta$ cohomologous to $\omega$. Then one
has:
$$\sigma(c\otimes c)=1,\quad \sigma(1\otimes h)=\sigma(h\otimes
1)=\varepsilon(h)$$
$$\sigma(x_j\otimes c)=0=\sigma(c\otimes x_j)$$
$$\sigma(x_i\otimes x_j)=\cases{m_{ij}&if $i<j$,\cr
m_{ji}& if $i>j$,\cr m_{jj}& if $i=j$,}$$ together with
relations depending on the recurrence relations for $\omega$ in
\cite{PVO2} and the value of $\theta$. Since $\theta\in Z(E(n)^*)$,
$\sigma$ is again a lazy cocycle (see \cite{lazy}). The matrix $L$ whose
$(i,j)$ entry is $\sigma(x_i\otimes x_j)$ is symmetric. Conversely, for any
symmetric matrix $L=(l_{ij})$ we will denote by $\sigma(L)$ the lazy
cocycle with $\sigma(x_i\otimes x_j)=l_{ij}$ contructed as above. The
corresponding cleft extension is the generalized Clifford algebra
with generators $u$ and $v_i$, for $i=1,\,\ldots n$, and relations
\begin{eqnarray}\label{privileged}
& u^2=1,\ uv_i+v_iu=0, \ & v_j^2=l_{jj}, \
v_iv_j+v_jv_i=2l_{ij} \quad \mbox{ for $i\neq j$}\quad
\end{eqnarray}
and with comodule algebra structure given by
\begin{equation}\label{privileged2}
\rho(u)=u\otimes c,\quad \rho(v_j)=1\otimes x_j+v_j\otimes c.
\end{equation}

For a symmetric matrix $L=(l_{ij})$ we denote by $Cl(L)$ the
algebra with generators and relations as in (\ref{privileged})
associated to the lazy cocycle $\sigma(L)$.

\section{$E(n)$-actions on central simple algebras}

In this section we analyze the action of $E(n)$ on central simple
algebras.

\begin{proposition}\label{azumaya}
Let $A$ be a central simple algebra which is an $E(n)$-module
algebra with action $\lu$. Then, there exist elements
$u,\,w_1,\,\ldots,\,w_n\in A$ such that $c \lu a=uau^{-1}$ and
$x_i \lu a=w_i(c\lu a)-aw_i$ for $i=1,\,\ldots n,$ and for all $a
\in A$. These elements satisfy the relations
$$u^2=\alpha,\quad w_i^2=l_{ii}\quad w_iu+uw_i=2\mu_i, \quad
w_iw_j+w_jw_i=2l_{ij}$$ for certain $\alpha\in \km$, and
$\mu_i,l_{ij} \in k$ for $1\le i\le j\le n.$
\end{proposition}
\pf By the Skolem-Noether Theorem for Hopf algebras, the action
of $E(n)$ on $A$ is inner. Hence there exists a convolution
invertible map $\pi\colon E(n)\to A$ such that $h\lu a=\sum
\pi(\h1)a\pi^{-1}(\h2)$ for all $h\in H,a \in A$. Putting
$u:=\pi(c)$ and $w_i:=-\pi^{-1}(x_i)=\pi(x_i)u^{-1}$ we get
$$c\lu a =uau^{-1},\qquad x_i\lu a=w_iuau^{-1}-aw_i,$$
for all $a \in A$. Since the action of $c^2$ is trivial, $u^2$ is
central in $A$ and therefore it belongs to $k$. Moreover, $u^2=\alpha$ is
nonzero because $u$ is invertible, so $\alpha\in\km$. \par
\smallskip
From $(x_ic+cx_i)\lu a=0$ for all $a \in A$, it follows that
$u^{-1}w_i+w_iu^{-1}$ is central in $A$. Thus, there is $\mu_i
\in k$ such that $w_iu+uw_i=2\mu_i$ for all $i=1,...,n$. Since
$(x_ix_j+x_jx_i)\lu a=0$ for every $a\in A$, one obtains that
$w_iw_j+w_jw_i$ belongs to the center of $A$. Then for every pair
$i,\,j$ there is $l_{ij}\in k$ such that $w_iw_j+w_jw_i=2l_{ij}$.
Finally, since $x_i^2$ acts as $0$ on $A$, we get $w_i^2=l_{ii}
\in k$. \qed

\begin{corollary}\label{clifford}
Let $A$ be a central simple algebra which is an $E(n)$-module
algebra. Then, there exists a submodule algebra $A'$ of $A$ which
is a quotient of a generalized Clifford algebra.
\end{corollary}
\pf The subalgebra $A'$ of $A$  generated by $u$ and the $w_j$'s
of Proposition \ref{azumaya} is a quotient of the generalized
Clifford algebra $Cl(\alpha,\mu,\,L)$ costructed from
the basis $\{v_0,\,\ldots,\,v_n\}$ of the vector space $V$ with
bilinear form associated to the matrix
\begin{equation}\label{matrice}
M= \left(\begin{array}{c|c}
\alpha&\mu_1~\mu_2~\cdots~\mu_n\\
\hline
\begin{array}{c}
\mu_1\\
\mu_2\\
\cdots\\
\mu_n
\end{array}&L
\end{array}
\right).
\end{equation}
\qed

\begin{lemma}\label{changes}
With notation as in Proposition \ref{azumaya}, the subalgebra
generated by $u$ and the $w_j$'s does not depend on the choice of
$\pi$. There exists a family $\pi_t\colon E(n)\to A$ with $t\in
\km$ for which $\mu_j=0$ for every $j$. For all values of $t$ the
matrix $L$ as in (\ref{matrice}) is uniquely determined.
\end{lemma}
\pf The action of $x_j$ and $c$ does not depend on the choice of
the map $\pi$. If we choose another map $\pi'\colon E(n)\to A$
inducing the same action, centrality of $A$ forces
$u'=\pi'(c)=tu$ for some $t\in\km$ and $w_j':=-\pi'(x_j)=w_j+s_ju$
for scalars $s_j$ for every
$j=1,\,\ldots,\,n$. Conversely, for any $t \in \km$ and
$s_j \in k$ for $j=1,\,\ldots,\,n$ the map $\pi'\colon E(n)\to A$
given by 
$\pi'(c)=tu,\; \pi^{-1}(x_j)=-(w_j+s_ju)$ and $\pi'(c^ax_{i_1}\cdots
x_{i_s})=\pi'(u)^a\pi'(x_{i_1})\cdots \pi'(x_{i_s})$ is
well-defined and satisfies $h\lu a=\sum \pi'(\h1)a\pi'^{-1}(\h2)$. Since
$$(w_j+s_ju)tu+tu(w_j+s_ju)=2(\mu_j+s_j\alpha)$$ the only choice for
$-(\pi')^{-1}(x_j)$ that skew-commutes with $\pi'(c)$ is for
$s_j=-\mu_j \alpha^{-1}$, with freedom on $t$. This choice is
always possible because $\alpha$ is invertible.  Hence, one has a
family $\pi_t\colon E(n)\to A$ given by $\pi_t(c)=tu$ and
$\pi_t(x_j)=w_ju+\mu_j\alpha^{-1}u$ and $\pi_t(c^ax_{i_1}\cdots
x_{i_s})=\pi_t(u)^a\pi_t(x_{i_1})\cdots \pi_t(x_{i_s})$. The
relations among the new $w_j$'s do not depend on $t$. \qed

Let us observe that if $u$ and $w_j$ are as in Lemma
\ref{changes} the action of $E(n)$ on the subalgebra generated by
$u$ and the $w_j$'s is as follows:
$$\begin{array}{lll}
c\lu u=u, & \quad x_j \lu u=-2uw_j; & \\
c\lu w_j=-w_j, & \quad x_j\lu w_i=-2l_{ij}.
\end{array}$$
In particular, the subalgebra generated by $u$ and the $w_j$'s is
a submodule algebra and it is independent of the choice of $\pi$.
We shall call it the subalgebra of $A$ {\em induced by the
$E(n)$-action} and we shall denote it by $Ind(A)$. The scalar
$\alpha$ and the symmetric matrix $L$ are called the {\it
invariants associated to} $A$. Lemma \ref{changes} establishes
that the matrix $M$ in (\ref{matrice})
can always be chosen to be $\left(\begin{array}{cc} \alpha&0\\
0&L \end{array}\right)$.

\begin{lemma}\label{strongly}
Let $A$ be a central simple algebra which is an $E(n)$-module
algebra. The action of $E(n)$ on $A$ is strongly inner if and only
if $\alpha$ is a square in $k$ and the matrix $M$ as in
(\ref{matrice}) has rank $1$.
\end{lemma}
\pf The action of $E(n)$ on $A$ is strongly inner if and only if
there exists a choice of $\pi$ which is an algebra map. Such a
$\pi$ has to be found among the $\pi_t$'s in Lemma \ref{changes}.
Suppose that for a suitable $t\in k^*$ the map $\pi_t$ is an
algebra map. Then for this $t$ one has:
$$\begin{array}{l}
(tu)^2=t^2\alpha=1; \qquad\quad  (w_j+s_ju)^2=l_{jj}-\mu_j^2\alpha^{-1}=0;\vspace{3pt} \\
(w_i+s_iu)(w_j+s_ju)+(w_j+s_ju)(w_i+s_iu)=2(l_{ij}+s_j\mu_i+
s_i\mu_j+s_is_j)=0.
\end{array}$$
Hence $\alpha$ must be a square, $l_{jj}\alpha=\mu_j^2$ and
$l_{ij}\alpha=\mu_i\mu_j.$ Thus the action of $E(n)$ on $A$
is strongly inner if and only if $\alpha$ is a nonzero square and
all minors of order $2$ of $M$ involving the first row and the
first column are $0$. Since $M$ is symmetric and
$\alpha=m_{11}\not=0$ it follows  that the action is strongly
inner if and only if $rk(M)=1$. \qed

Let us observe that Lemma \ref{changes} and Lemma \ref{strongly}
imply that in this case we can always make sure that the only
nonzero entry in the matrix $M$ in (\ref{matrice}) is
$m_{11}=\alpha=1$, i.e., that $L=0$.

\begin{lemma}\label{bq}
Let $H$ be a Hopf algebra and let  $A$ and $B$ be Yetter-Drinfeld
$H$-module algebras. If $H$ acts on $B$ in a strongly inner way,
then $A\# B\cong A\otimes B$ as algebras. In particular, if $A$
and $A'$ are Brauer equivalent $H$-Azumaya algebras and $A$ is
central simple then $A'$ is also central simple.
\end{lemma}
\pf Let $f\colon H\to B$ be a convolution invertible algebra map
such that $h\lu b=\sum f(\h1)bf^{-1}(\h2)$ for all $h\in H, b\in
B$. We check that the linear map
$$\Lambda\colon A\otimes B\to A\# B,\ a\otimes b \mapsto \sum\a0\# f(\a1)b$$
is an algebra map. For $a,\,c\in A$ and $b,\,d\in B,$ we have:
$$\begin{array}{ll}
\Lambda(a\otimes b)\Lambda(c\otimes d) & = \left(\sum\a0\#f(\a1)b
  \right)\left(\sum\c0\# f(\c1)d\right) \vspace{3pt}\\
 & = \sum\a0\c0\#[\c1\lu\left(f(\a1)b\right)]f(\c2)d \vspace{3pt}\\
 & = \sum\a0\c0\# f(\c1)\left(f(\a1)b\right) f^{-1}(\c2)f(\c3)d \vspace{3pt}\\
 & = \sum\a0\c0\# f(\c1)f(\a1)bd \vspace{3pt} \\
 & = \sum\a0\c0\# f(\c1\a1)bd \vspace{3pt}\\
 & = \Lambda((a\otimes b)(c\otimes d)).\\
\end{array}$$
The inverse of $\Lambda$ is the map $a\# b\to \sum\a0\otimes
f^{-1}(\a1)b$. \par \smallskip

If $A$ and $A'$ are as in the hypothesis, then $A\# End(P)\cong
A'\# End(Q)$ for some Yetter-Drinfeld $H$-modules $P$ and $Q$.
Since the action of $E(n)$ on $End(P)$ and $End(Q)$ is strongly
inner we have that $A\otimes End(P) \cong A'\otimes End(Q).$ Then
$A' \otimes End(Q)$ is central simple and by the Double
Centralizer Theorem, $A'=C(End(Q))$ is central simple. \qed

We shall use the above Lemma for the investigation of the Brauer
group $BM(k,\,E(n),\,R_A)$ in the next section.

\begin{proposition}\label{invariance}
Let $M \in M_n(k)$ and let $A$ be a central simple algebra which
is an $(E(n), R_M)$-module algebra. Let $u,\,w_j$ be the
generators of $Ind(A)$ with associated invariants $\alpha$ and
$L$. Let $B$ be a central simple algebra on which $E(n)$ acts in a
strongly inner way, endowed with the coaction induced by $R_M$.
Then the associated invariants for $Ind(A\#B)$ are again $\alpha$
and $L$.
\end{proposition}
\pf By Lemma \ref{bq}, $A\# B$ is also a central simple algebra,
so the action of $E(n)$ on $A\# B$ is again inner. Let $u,\,w_j$
be the generators for $Ind(A)$ and $v,\,v_j$ be the generators
for $Ind(B)$. In particular, $x_j\lu v_i=0$ for every $i$ and
$j$. It is not difficult to check that the map $\Lambda$ in Lemma
\ref{bq} preserves the $c$-action. Since
$$\Lambda^{-1}(c\lu a\# c\lu b)=(u\otimes v)\Lambda^{-1}(a\#
b)(u\otimes v)^{-1},$$ it follows that for $U:=\Lambda(u\otimes
v)$ we have:
$$c\lu(a\# b)=U(a\# b)U^{-1},\quad U^2=u^2\# v^2=\alpha.$$

Next we consider the braiding for the elements $w_j$. Let
$\Psi_0$ denote the braiding of two $E(n)$-modules induced by the
quasi-triangular structure $R_0$. For $b\in B$ we have:

\begin{quote}
$\begin{array}{l}
\hspace{-1cm} \Psi_M(w_j\otimes b)= \sum R_M^{(2)}\lu b\otimes R_M^{(1)}\lu w_j \vspace{3pt} \\
\hspace{0.7cm} = \Psi_0(w_j\otimes b)+
{1\over2}\sum_{P}(-1)^{{|P|(|P|-1)}\over2}\sum_{F,\; |F|=|P|,\;
\eta\in S_{|P|}}sign(\eta)m_{P,\eta(F)} \vspace{3pt}\\
\hspace{1.2cm}  \bigl[(c^{|P|}x_F\lu b)\otimes (x_P\lu w_j)+(c^{|P|}x_F\lu b) \otimes (cx_{P}\lu w_j)+ \vspace{3pt} \\
\hspace{1.2cm} (c^{|P|+1}x_F\lu b) \otimes (x_P\lu
w_j)-(c^{|P|+1}x_F\lu b)\otimes (cx_{P}\lu w_j\bigr)].
\end{array}$
\end{quote}
By the particular module structure of $Ind(A)$ the terms of the
form $cx_{P}\lu w_j$ with $|P|>1$ will vanish. Hence
\begin{quote}
$\begin{array}{ll} \Psi_M(w_j\otimes b) & = \Psi_0(w_j\otimes b)+
{1\over2}\sum_{p,f}m_{pf} \bigl[(cx_f\lu b) \otimes (x_p\lu
w_j) \vspace{3pt} \\
 &\ \   +(cx_f\lu b )\otimes (cx_p\lu w_j)+ (x_f\lu b) \otimes (x_p\lu w_j) \vspace{3pt} \\
 &\ \   -(x_f\lu b) \otimes (cx_p\lu w_j) \bigr] \vspace{3pt} \\
 & =\Psi_0(w_j\otimes b)-2\sum_{p,f}m_{pf}l_{pj}(cx_f\lu b)\otimes 1.
\end{array}$
\end{quote}
Analogously,
$$\Psi_M(b\otimes w_j)= \Psi_0(b\otimes w_j)-2\sum_{p,f}m_{pf}l_{fj}(1\otimes (x_p\lu
b)).$$ A similar computation shows that $\Psi_M(b\otimes
v_j)=\Psi_0(b\otimes v_j)$ and $\Psi_M(v_j\otimes
b)=\Psi_0(v_j\otimes b).$ Let us consider the elements
$$W_i:=w_i\#1+1\#v_i-2\sum_{p,f}m_{pf}l_{fi}(1\# v_p)$$
for $i=1,\,\ldots,\,n$. Then, for homogeneous $a$ and $b$ with
respect to the ${\mathbb Z}_2$-grading induced by the action of
$c$ we have:
\begin{quote}
$\begin{array}{l}
\hspace{-1cm} W_i(c\lu(a\#b))-(a\#b)W_i = \vspace{3pt} \\
\hspace{1cm} = (w_i\#1)(c\lu a\#c\lu b)+(1\#v_i)((-1)^{\deg(a)}a\#c\lu b)\vspace{3pt} \\
\hspace{1.5cm} -2\sum_{p,f}m_{pf}l_{fi}(1\#v_p)((-1)^{\deg(a)}a\# c\lu b))-(a\#b)(w_i\#1)\vspace{3pt} \\
\hspace{1.5cm} -a\# bv_i+2(a\#b)\sum_{p,f}m_{pf}l_{fi}(1\#v_p) \vspace{3pt} \\
\hspace{1cm} = w_i(c\lu a)\#(c\lu b)-(-1)^{\deg(b)}aw_i\# b+a\#v_i(c\lu b)-a\#bv_i\vspace{3pt}\\
\hspace{1.5cm} +2\sum_{p,f}m_{pf}l_{fi}a\#[v_p(c\lu b)-bv_p-x_p\lu b] \vspace{3pt} \\
\hspace{1cm} = (x_i\lu a)\#(c\lu b)+a\# (x_i\lu b)=x_i\lu(a\#b).
\end{array}$
\end{quote}
Hence $Ind(A\# B)$ is the subalgebra 
generated by $U$ and the $W_i$'s. Besides, $UW_i+W_iU=((c\lu
W_i)+W_i)U=0$. It may be checked that
$$W_i^2=-\frac{1}{2}(x_i\lu W_i)=l_{ii}\quad {\normalfont and} \quad
W_iW_j+W_jW_i=-\frac{1}{2}(x_i\lu
W_j)=l_{ij}.~\hfill\Box$$\par\bigskip

\begin{corollary}\label{br-invariance}
With notation as in Proposition \ref{invariance} the scalar
$\alpha$ and the matrix $L$ are invariant under the Brauer
equivalence induced by $R_M$.
\end{corollary}
\pf Apply Proposition \ref{invariance} with $B=End(P)$ for a
$E(n)$-module $P$.~\qed

The following lemma will be needed in the forthcoming section.

\begin{lemma}\label{surjective}
Let $\sigma$ be a $2$-cocycle of a Hopf algebra $H$ and let $f$
be the map $f\colon H\to End(H)$ defined by
$$(f(h))(a)=\sum\sigma(\h1\otimes\a1)\h2\a2$$ for every $a,\,h\in H$. Then:
\begin{enumerate}
\item[(i)] The map $f$ is convolution invertible with inverse
$$f^{-1}(h)=\sum \sigma^{-1}(S(\h2)\otimes\h3)f(S(\h1))$$ for every
$h\in H$;
\item[(ii)] For every $h,l\in H$ there holds $f(h)\circ
f(l)=\sum\sigma(\h1\otimes\l1)f(\h2\l2)$; 
\item[(iii)] $H$ measures $End(H)$ by means of $h\lu F=\sum f(\h1)\circ
F\circ f^{-1}(\h2)$ for every $F\in End(H)$.
\end{enumerate}
\noindent Let now $\sigma$ be a lazy cocycle. Then:
\begin{enumerate} 
\item[(iv)] The weak action $\lu$ is indeed an action;
\item[(v)] The induced subalgebra on $End(H)$ is a quotient of the cleft extension of $k$
determined by $\sigma$.
\end{enumerate}
\end{lemma}
\pf {\it (i)} This assertion follows from \cite[Proposition
7.2.7]{Mon} which holds for all crossed products. Indeed $f(h)$
is just multiplication by Blattner and Montgomery's $\gamma(h)$
in $_{\sigma}H=k\#_\sigma H$ in their notation. Since $\gamma\colon
H\to k\#_\sigma H$ given by $h\to 1\#h$ is convolution invertible
with convolution inverse $\mu(h)=\sum
\sigma^{-1}(S(\h2)\otimes\h3)S(\h1)$, the convolution inverse for
$f$ is $g\colon H\to End(H)$ given by $g(h)=\sum
\sigma^{-1}(S(\h2)\otimes\h3)f(S(\h1))$ for every $h \in H$. \par
\smallskip

{\it (ii)} The second assertion is also essentially contained in
\cite[Proposition 7.2.7]{Mon} because of the description of
$f(h)$ as multiplication by $\gamma(h)$. We can simply view $f$
as the left regular representation of $_\sigma H$. \par \smallskip

{\it (iii)} By construction $\lu$ defines a measuring. \par
\smallskip

{\it (iv)} It follows from {\it (i)} and {\it (ii)} applied to $f(S(h))$ that

$$\begin{array}{ll}
f^{-1}(hl) & = \sum \sigma^{-1}(S(\l2)S(\h2)\otimes\h3\l3)f\left(S(\l1)S(\h1)\right) \vspace{3pt} \\
 & =\sum \sigma^{-1}(S(\l3)S(\h3)\otimes\h4\l4)\sigma^{-1}(S(\l2)\otimes S(\h2)) \vspace{3pt}\\
 &\quad  f(S(\l1))\circ f(S(\h1)). \vspace{3pt}
\end{array}
$$
Since $\sigma^{-1}$ is a right cocycle
$$
\begin{array}{ll}
f^{-1}(hl) & =\sum\sigma^{-1}(S(\h2)\otimes \h5\l4)\sigma^{-1}(S(\l2)\otimes
  S(\h3)\h4\l3) \vspace{3pt} \\
& \quad f(S(\l1))\circ f(S(\h1))\vspace{3pt}\\
& = \sum\sigma^{-1}(S(\h2)\otimes \h3\l2)f^{-1}(\l1)\circ f(S(\h1)).\\
\end{array}$$
Using this equality, we obtain:
$$\begin{array}{l}
\sum f^{-1}(\h1\l1)\sigma^{-1}(\h2\otimes\l2)=\vspace{3pt}\\
\hspace{0.7cm}=\sum\sigma^{-1}(S(\h2)\otimes \h3\l2)\sigma^{-1}(\h4\otimes\l3)f^{-1}(\l1)\circ f(S(\h1))\vspace{3pt}\\
\hspace{0.7cm}=\sum \sigma^{-1}(S(\h2)\otimes\h5)\sigma^{-1}(S(\h3)\h4\otimes\l2)f^{-1}(\l1)\circ f(S(\h1))\vspace{3pt}\\
\hspace{0.7cm} =f^{-1}(l)\circ f^{-1}(h).
\end{array}$$
Therefore,  $$f^{-1}(hl)=
\sum\sigma(\h2\otimes\l2)f^{-1}(\l1)\circ f^{-1}(\h1).$$ If
$\sigma$ is lazy we also have:
$$f(hl)=\sum \sigma^{-1}(\h2\otimes\l2)f(\h1)\circ
f(\l1),$$ $$ f^{-1}(hl)= \sum\sigma(\h1\otimes\l1)f^{-1}(\l2)\circ
f^{-1}(\h2).$$

Let $\lu$ be the weak action defined above. Then $1\lu F=F$ for
every $F\in End(H)$. Besides, for $F\in End(H)$ and $h,\,l\in H$
we have:
$$\begin{array}{ll}
(hl)\lu F & = \sum f(\h1\l1)\circ F\circ f^{-1}(\h2\l2) \vspace{3pt} \\
& =\sum f(\h1)\circ f(\l1)\sigma^{-1}(\h2\otimes\l2)\circ F\circ
\sigma(\h3\otimes\l3) \vspace{3pt}\\
& \quad f^{-1}(\l4)\circ f^{-1}(\h4)\vspace{3pt}\\
& =\sum (f(\h1)\circ (f(\l1)\circ F\circ f^{-1}(\l2))\circ f^{-1}(\h2))\vspace{3pt}\\
& =h\lu(l\lu F).
\end{array}$$
Hence $End(H)$ is an $H$-module algebra. \par \smallskip

{\it (v)} Since $End(H)$ is a central simple algebra, the induced
subalgebra, i.e., the subalgebra of $End(H)$ generated by the
elements of the form $f(h)$ for $h\in H$ is well-defined and by
the relations for $f(hl)$ we deduce that it is a quotient of
$H_\sigma\cong _{\sigma}\!H$, the module algebra twist of $H$ by
$\sigma$. \qed

Given a lazy cocycle $\sigma$ we shall denote the $H$-module
algebra structure on  $End(H)$ above defined by $A^\sigma$. 
%

\section{Subgroups of $BQ(k,\,E(n))$}

\subsection{Subgroups arising from (co)quasi-triangular structures}

We first focus on the computation of the Brauer group
$BM(k,\,E(n),R)$. Since $E(n)^{op}\cong E(n),$ by Corollary
\ref{classi} and Remark \ref{dual}, we can reduce the computation
of $BM(k,\,E(n), R_M)$ for $M\in M_n(k)$ to the computation of
$BM(k,\,E(n), R_{J_l})$ with $J_l$ as in $($\ref{matrici}$)$. In
particular we may always assume that the matrix $M$ is
skew-symmetric, as we will do in the sequel.

\begin{lemma}\label{gruppo}
Let $n$ and $r$ be nonnegative integers with $r \le n$. Let $M
\in M_n(k)$ be skew-symmetric and such that the last $r$ rows and
columns are $0$. Let $\sym$ be the set of symmetric matrices of
the form
$$L=\begin{array}{c}
\begin{array}{cc}
\overbrace{}^{n-r}&\overbrace{}^r
\end{array}\\
\left(\begin{array}{c|c}
0&L_1\\
\hline L_1^t&L_2
\end{array}\right).
\end{array}
$$ Then the operation
$$L\oplus N=L+N-2(NML)-2(NML)^t=L+N-2(NML)+2(LMN)$$ turns
$\sym$ into a group with unit $0$ and with opposite $-L$.
\end{lemma}
\pf It is clear that $0$ is the neutral element with respect to
this operation and that $-L$ is a right and left opposite for
$L$. To check associativity, let $L$, $N$ and $P\in \sym$. Then,
$$\begin{array}{l}
(L\oplus N)\oplus P =(L+N-2(NML)+2(LMN))\oplus P \vspace{3pt} \\
\hspace{1.5cm} = L+N-2(NML)+2(LMN)+P-2PM(L+N-2(NML) \vspace{3pt} \\
\hspace{1.8cm} +2(LMN))+2(L+N-2(NML)+2(LMN))MP \vspace{3pt} \\
\hspace{1.5cm} = L+N+P-2(NML)+2(LMN)-2(PML) \vspace{3pt} \\
\hspace{1.8cm} -2(PMN)+2(LMP)+2(NMP) \vspace{3pt} \\
\hspace{1.5cm} = L+N+P-2(PMN)+2(NMP)-2((N+P-2(PMN)\vspace{3pt} \\
\hspace{1.8cm} +2(NMP))ML)+2(LM(N+P-2(PMN)+2(NMP)) \vspace{3pt} \\
\hspace{1.5cm} = L\oplus (N+P-2(PMN)+2(NMP)) \vspace{3pt} \\
\hspace{1.5cm} = L\oplus(N\oplus P).
\end{array}$$
We have used that $MTM=0$ for every $T\in\sym$. This holds because
of the particular block form of $M$. Hence $(\sym,\oplus)$ is a
group. \qed

\begin{lemma}\label{central}
With hypothesis as in Lemma \ref{gruppo}, the group $Sym_{M,n,r}(k)$
is a central extension with kernel isomorphic to $(Sym_r(k), +)$
and quotient isomorphic to $(M_{n-r,n}(k), +)$.
\end{lemma}
\pf Let $M'$ denote the submatrix of $M$ corresponding to the first
$(n-r)$ rows and columns. Writing every element $L=\left(\begin{array}{cc}
0&L_1\\
L_1^t&L_2
\end{array}\right)\in Sym_{M,n,r}(k)$ as a pair $(L_1,\,L_2)\in
M_{n-r,r}(k)\times Sym_r(k)$, the product in $Sym_{M,n,r}(k)$
becomes:
$$(L_1,L_2)\oplus(N_1,N_2)=(L_1+N_1,L_2+N_2-2N_1^tM'L_1+2L_1^tM'N_1)$$
so the elements of type $(0, S)$ form a subgroup isomorphic to
$(Sym_r(k),+)$. One can directly see that this subgroup is
central and that the image of the projection $(L_1,L_2)\mapsto L_1$
is isomorphic to $(M_{n-r,n}(k), +)$. \qed

\begin{theorem}\label{reduction}
Let $n\ge1,$ let $r\le n$ and $M$ be an $n\times n$
skew-symmetric matrix whose last $r$ rows and $r$ columns are
zero. Let $M'\in M_{n-r}(k)$ be the submatrix 
of $M$ corresponding to the first
$(n-r)$ rows and columns. Then
there is a split short exact sequence:
$$1\longrightarrow \sym\longrightarrow
BM(k,\,E(n), R_M)\longrightarrow BM(k,\,E(n-r),
R_{M'})\longrightarrow 1.$$ If $n=r$, we have
$E(n-n)=E(0)=k{\mathbb Z}_2$, $M'=0$ and $BM(k,\,E(0), R_{M'})=BW(k)$,
the Brauer-Wall group of $k$.
\end{theorem}
\pf The natural injection $j_{n-r}\colon (E(n-r),R_{M'})\to(E(n),
R_M)$ is a quasi-triangular map, hence the braidings in the
corresponding categories of modules coincide, due to the
particular form of $M$ and $M'$. Indeed, if we look at the action
of $\sum R_M^{(2)}\otimes R_M^{(1)}$ on tensor products we see
that the action of the $x_j$'s with $j>r$ never occurs because
all the monomials involving these elements appear with zero
coefficient. Therefore the pull-back along $j_{n-r}$ induces a
group homomorphism $j_{n-r}^*\colon BM(k,\,E(n), R_M)\to
BM(k,\,E(n-r), R_{M'})$. Since the natural projection
$p_{n-r}\colon (E(n), R_{M})\to(E(n-r), R_{M'})$ is a
quasi-triangular map, the pull-back along $p_{n-r}$ induces a
group homomorphism $p_{n-r}^*\colon BM(k,\,E(n-r), R_{M'})\to
BM(k,\,E(n), R_{M})$ splitting $j_{n-r}^*$. Let us compute
$Ker(j_{n-r}^*)$. Its elements are represented by endomorphism
algebras with strongly inner $E(n-r)$-action. Let $A=End(P)$ be
such a representative. The matrix $L$ associated to $A$ is
uniquely determined by Lemma \ref{changes} and by Corollary
\ref{br-invariance} and it is invariant under Brauer equivalence
induced either by $R_M$ and by $R_{M'}$. By Lemma \ref{strongly}
applied to $E(n-r)$ there holds  $l_{ij}=0$ for $1\le i,\,j\le r$, i.e.,
$L$ lies in $\sym$ and we can make sure that the associated
scalar $\alpha=1$. Hence we have a well-defined map $$\chi\colon
Ker(j_{n-r}^*)\to \sym$$
$$[A]\mapsto L.$$
We next check that $\chi$ is a group homomorphism. For this
purpose let $B$ be a representative of another element in
$Ker(j_{n-r}^*)$. Assume that $Ind(B)$ is generated by $v$ and
$v_i$ as in Lemma \ref{changes}, and let $L'=(l'_{ij})$ be a
matrix such that $\chi([B])=L'$. Then
$$\Psi_M(v_j\otimes a) =\Psi_0(v_j\otimes a)-2\sum_{p,f\le
n-r}m_{pf}l'_{pj}(cx_f\lu a)\otimes 1,$$ for all $a\in A$.
Similarly,
$$\Psi_M(a\otimes v_j)= \Psi_0(a\otimes v_j)-2\sum_{p,f\le
n-r}m_{pf}l'_{fj}(1\otimes (x_p\lu a)),$$ while the braiding
for the $w_j$'s is as in Proposition \ref{invariance}. Again, by
Proposition \ref{invariance}, the element $U=\Lambda(u\otimes v)$
satisfies: $U^2=1$ and $c\lu(a\#b)=U(a\#b)U^{-1}$. \par \smallskip

Let
$$W_j:=w_j\#1+1\#v_j-2\sum_{p,f}m_{pf}l_{fj}\#v_p-2
\sum_{p,f}m_{pf}l'_{pj}w_f\#1.$$ We consider the $\Ent_2$-grading on
$A,B$ induced by the $c$-action. For homogeneous elements $a\in
A, b\in B$ one has: 

$$\begin{array}{l}
W_j(c\lu(a\#b))-(a\#b)W_j= \vspace{3pt} \\
\hspace{0.3cm} =w_j(c\lu a)\#(c\lu b)+a\#v_j(c\lu b) \vspace{3pt} \\
\hspace{0.6cm} +2\sum_{p,f}m_{pf}l'_{pj}((x_f\lu a)\#(c\lu
b))-2\sum_{p,f}m_{pf}l_{fj}(a\#v_p(c\lu b)) \vspace{3pt} \\
\hspace{0.3cm}
-4\sum_{p,f}m_{pf}l_{fj}\sum_{r,s}m_{rs}l'_{rp}((x_s\lu
a)\#(c\lu b)) \vspace{3pt} \\
\hspace{0.6cm} -2\sum_{p,f}m_{pf}l'_{pj}(w_f(c\lu a)\#(c\lu
b))-aw_j\#(c\lu b) \vspace{3pt} \\
\hspace{0.6cm}+2\sum_{p,f}m_{pf}l_{fj}(a\# (x_p\lu b))-a\# bv_j \vspace{3pt} \\
\hspace{0.6cm}+2\sum_{p,f}m_{pf}l_{fj}(a\#bv_p)+2\sum_{p,f}m_{pf}l'_{pj}
(aw_f\#(c\lu b))\vspace{3pt} \\
\hspace{0.6cm}-4\sum_{p,f}m_{pf}l'_{pj}\sum_{r,s}m_{rs}l_{sf}(a\#
(x_r\lu b)) \vspace{3pt} \\
\hspace{0.3cm} =((x_j\lu a)\# (c\lu b)+a\#(x_j\lu b))\vspace{3pt}\\
\hspace{0.6cm}+ 2\sum_{p,f}m_{pf}l'_{pj}[(x_f\lu a-w_f(c\lu
a)+aw_f)\#(c\lu b)]\vspace{3pt} \\
\hspace{0.6cm}+2\sum_{p,f}m_{pf}l_{fj}[a\#(x_p\lu b-v_p(c\lu b)+bv_p)]\vspace{3pt} \\
\hspace{0.6cm}-4\sum_s(LM^tL'M)_{js}((x_s\lu a)\# (c\lu b))
-4\sum_s(L'MLM^t)_{js}(a\#(x_s\lu b))\vspace{3pt} \\
\hspace{0.3cm} =x_j\lu(a\#b)
\end{array}$$
because $MTM$ for every $T\in\sym$.
Therefore $Ind(A\# B)$ is the subalgebra generated by $U$ and the
$W_j$'s. Since $c\lu W_j=-W_j$ we have $UW_j+W_jU=0$. Besides,
$$\begin{array}{ll}
W_j^2 & =-\frac{1}{2}(x_j\lu W_j) \vspace{3pt}\\
 & =l_{jj}+l'_{jj}-2\sum_{p,f}m_{pf}l_{fj}l'_{pj}-2
 \sum_{p,f}m_{pf}l'_{pj}l_{fj}\vspace{3pt} \\
& =(l_{jj}+l'_{jj})-2(L'ML)_{jj}-2(L'ML)_{jj}.
\end{array}$$
Similarly,
$$(W_jW_i+W_iW_j)=-\frac{1}{2}(x_j\lu W_i)
=(l_{ij}+l'_{ij})-2(L'ML)_{ji}-2(L'ML)_{ij}.$$ Hence $\chi\colon
Ker(j_{n-r}^*)\to\sym$ is a group homomorphism. We claim that it
is injective. If $\chi([A])=0$ then the action of $E(n)$ on the
endomorphism algebra $A$ is strongly inner, forcing $[A]=[k]$. We
finally show that $\chi$ is surjective. To this end let $L\in
\sym$ and let $\sigma$ be a lazy cocycle with $\sigma(x_i\otimes
x_j)=-l_{ij}$. Such a cocycle exists because of the results in
Section \ref{twistings}. Let $A^{\sigma}$ be constructed as in
Lemma \ref{surjective}. For the induced subalgebra
$Ind(A^\sigma)$ one has $f(h)\circ
f(l)=\sum\sigma(\h1\otimes\l1)f(\h1\l1)$. Then, for every
element of $E(n-r)$, we have $f(h)f(k)=f(hk)$, i.e, the action of
$E(n-r)$ is strongly inner. Therefore, $E(n)$, as a vector space,
is a $E(n-r)$-module (hence, via $R_{M'}$, a Yetter-Drinfeld
module). Thus $A^{\sigma}$ is a $(E(n-r),\,R_{M'})$-Azumaya
algebra. Since the braiding induced by $R_M$ does not involve the
action of the $x_j's$ for $j>n-r$, the algebra $A^{\sigma}$ is
also $(E(n),\,R_{M})$-Azumaya. Besides, for $u=f(c)$ we have
$u^2=1$; for $v_j=f(x_j)$ we have $v_j^2=-l_{jj}$, and
$v_iv_j+v_iv_j=-2l_{ij}$. Then for $w_j=v_ju$ for $j=1,\,\ldots,
n$ we have $w_j^2=l_{jj}$ and $w_iw_j+w_jw_i=2l_{ij}$. Hence
$\chi([A^{\sigma}])=L$ and $\chi$ is surjective. \qed

\begin{theorem}\label{direct}
Let $A\in Sym_n(k)$. Then $BM(k,\,E(n),\,R_A)$ is isomorphic to
the direct product of $(Sym_n(k),+)$ and $BW(k)$.
\end{theorem}
\pf If $A$ is symmetric, then $BM(k,E(n),\,R_A)\cong
BM(k,\,E(n),\,R_0)$ and it is an abelian group because the
category of $E(n)$-modules endowed with the braiding stemming
from $R_A$ is symmetric. Now apply Proposition \ref{reduction} for
$M=0$ and $r=n$. In this case, $\sym\cong Sym_n(k)$, hence the
proof. \qed

\begin{theorem}\label{semi}
Let $n\ge1$, $j\le\left[\frac{n}{2}\right]$ and let $J_{l}$ and
$J_{l}'$  be the square matrices of order $n$ and $2l$ respectively,
defined as in (\ref{matrici}). Then  $BM(k,E(n),\,R_{J_l})$ is a
semidirect product of $Sym_{n,n-2l,J_l}(k)$ and
$BM(k,\,E(2l),R_{J_l'})$.
\end{theorem}
\pf Apply Proposition \ref{reduction} to $M=J_l$ and $r=n-2l$.
\qed

We point out that $BM(k,E(n),R_{J_l})$ for $l>0$ is not abelian
because its subgroup $Sym_{n, n-2l, J_l}(k)$ is not. Theorem
\ref{semi} together with the classification of orbits in Section
\ref{twistings} shows that in order to compute the Brauer groups
$BM(k,\,E(n), R_A)$ for every matrix $A$ it is enough to
understand those of type $BM(k,\,E(2l), R_{J_l})$. This latter
group seems to be more complicated than those corresponding to
$R_0$. The group $BM(k,\,E(2l), R_{J_l})$ contains two copies of
$BM(k,E(l), R_0)$, obtained as the images of splitting maps
similar to $p_{n-r}^*$ in the proof of Proposition
\ref{reduction} isolating either $c$ and $x_1,\,\ldots,\,x_l$ or
$c$ and $x_{l+1},\,\ldots x_{2l}$. The intersection of this two
subgroups is the Brauer-Wall group of $k$ because it corresponds to
those $E(2l)$-Azumaya algebras with trivial action of all the
skew-primitives elements.

\subsection{Subgroups arising from the Hopf automorphisms group}

In \cite[Theorem 5]{VZ2} an exact sequence relating the Brauer
group of a Hopf algebra $H$ and its Hopf automorphism group
$Aut_{Hopf}(H)$ was constructed. We recall the construction of
this sequence and apply it to $E(n)$. Given $\alpha \in
Aut_{Hopf}(H)$ let $H_{\alpha}$ denote the right $H$-comodule $H$
with the following left $H$-action:
$$h\cdot m=\sum \alpha(h_{(2)})mS^{-1}(h_{(1)}) \qquad
\forall h \in H, m \in H_{\alpha}.$$ The algebra
$A_{\alpha}=End(H_{\alpha})$ is an $H$-Azumaya
algebra with $H$-module and $H$-comodule structures defined by

$$\begin{array}{ll}
(h\cdot f)(m)=\sum h_{(1)} \cdot f(S(h_{(2)})\cdot m), \vspace{3pt} \\
\rho(f)(m)=\sum f(m_{(0)})_{(0)} \otimes
S^{-1}(m_{(1)})f(m_{(0)})_{(1)}
\end{array}$$
for $h \in H, m \in H_{\alpha}, f \in A_{\alpha}.$ The map
$\pi:Aut_{Hopf}(H) \rightarrow BQ(k,H), \alpha \mapsto
[A_{\alpha^{-1}}]$ is a group homomorphism whose kernel is
$G(D(H))/G(D(H)^*)$, where $G(H)$ denotes the group of grouplike
elements of a Hopf algebra $H$. It is well-known that $G((D(H))=G(H)\times
G(H^*)$ and
$$\begin{array}{l}
G(D(H)^*)= \\
\quad \{(g,\lambda) \in G(H) \times G(H^*) \vert \sum gh_{(1)}
\lambda(h_{(2)})=\sum h_{(2)}g\lambda(h_{(1)})\ \forall h \in
H\}.
\end{array}$$
There is a group homomorphism $\theta:G(D(H)) \rightarrow
Aut_{Hopf}(H)$ defined by $\theta((g,\lambda))(h)=\sum
\lambda(h_{(1)})gh_{(2)}g^{-1}\lambda^{-1}(h_{(3)})$ for all $h
\in H$. It was proved in \cite[Theorem 5]{VZ2} that $\pi$ and
$\theta$ fit in the following exact sequence:
$$\xymatrix{1 \ar[r] & G(D(H)^*) \ar[r] & G(D(H)) \ar[r]^{\theta} &
Aut_{Hopf}(H) \ar[r]^{\pi} & BQ(k,H).}$$ We compute the first
three terms of this sequence for $H=E(n)$ and determine the
quotient of $Aut_{Hopf}(E(n))$ embedding in $BQ(k,E(n))$. The
case $E(1)$ was treated in \cite[Example 7]{VZ2}. Recall that the
group $Aut_{Hopf}(E(n))$ is isomorphic to $GL_n(k)$. Any $\alpha
\in Aut_{Hopf}(E(n))$ is of the form
$$\alpha(c)=c, \quad \alpha(x_i)=\sum_{j=1}^n t_{ij}x_j \quad \forall
i=1,\ldots ,n,$$ for $T_{\alpha}=(t_{ij}) \in GL_n(k)$. Since
$E(n)$ is self-dual, $G(E(n)^*)$ is easily computed. It consists
of the elements $\varepsilon$ and $C=1^*-c^*$. Then
$G(D(H))=\{(1,\varepsilon),(c,\varepsilon),(1,C),(c,C)\} \cong
\Ent_2 \times \Ent_2.$ It may be checked that
$G(D(H)^*)=\{(1,\varepsilon),(c,C)\} \cong \Ent_2$. The quotient
group $G(D(H))/G(D(H^*)) \cong \Ent_2$ is generated by the class
of $(c,\varepsilon)$ and $\theta(c,\varepsilon)$ is the
automorphism group corresponding the to matrix $-Id.$ Consider
$\Ent_2$ inside $GL_n(k)$ as the subgroup generated by $-Id$. Then
$GL_n(k)/\Ent_2$ embeds in $BQ(k,E(n))$. For $T \in GL_n(k)$ let
$\overline{T}$ denote the class of $T$ in $GL_n(k)/\Ent_2$. The
embedding is given by the map $\pi:GL_n(k)/\Ent_2 \rightarrow
BQ(k,E(n)),\overline{T} \mapsto [A_{T^{-1}}]$. \par \medskip

The group $Aut_{Hopf}(H)$ acts on $BQ(k,H)$ by setting $
[B]^{\alpha}=[A_{\alpha}][B][A_{\alpha^{-1}}]$ for any $\alpha \in
Aut_{Hopf}(H), [B] \in BQ(k,H)$. In \cite[Theorem 4.11]{CVZ2}, an
alternative description of this action is given:
$[B]^{\alpha}=[B(\alpha)]$ where $B(\alpha)=B$ as a $k$-algebra
but with $H$-module and $H$-comodule structures given by
$$h \cdot_{\alpha} b=\alpha(h) \cdot b,\qquad \rho(b)=\sum b_{(0)}
\otimes \alpha^{-1}(b_{(1)}),$$ for all $h \in H,b \in B(\alpha).$
We study this action for $H=E(n)$ and show that it stabilizes the
subgroup $Sym_n(k)$ found in the previous section. We will prove
that this action corresponds to the action of $GL_n(k)$ on
$Sym_n(k)$ given by $L^T=TLT^t$ and we will next embed the
subgroup $Sym_n(k) \rtimes (GL_n(k)/\Ent_2)$ into $BQ(k,E(n))$.
\par \smallskip

Let $L\in Sym_n(k)$  
and let $A^L:=A^{\sigma}$ as in the proof of Theorem \ref{reduction}.
Then for $f \in A^{L}$,
$$c \cdot f=U \circ f \circ U,\qquad x_i \cdot f=W_i \circ U \circ f \circ
U-f \circ W_i$$ where $U$ and $W_i$ are the linear maps
corresponding to the multiplication in $Cl(L)$ by $u$ and $w_i$
respectively. For $T=(t_{ij}) \in GL_n(k)$ we show that $A^L(T)$ is
again central and simple and we compute its induced
subalgebra by the $E(n)$-action. Since all
automorphisms of $E(n)$ fix $c$, the $E(n)$-comodule structure on
$A^{L}(T)$ is again determined by $R_0$ so that $A^{L}(T)$ is again a
representative of an element in $BM(k, E(n), R_0)$. Besides, the
$c$-action on $A^{L}(T)$ coincides with the $c$-action on $A^{L}$ so
the $c$-action is again strongly inner and $A^{L}(T)$ is again the
representative of an element in $Ker(j^*_0)=Sym_n(k)$. The
element induced by the $c$-action is again $U$ while the element
$W'_i$ induced by the $x_i$-action is $W'_i=\sum_{j=1}^n
t_{ij}W_j$. One may check that
$$(W'_i)^2=l'_{ii},\qquad W'_iW'_j+W'_jW'_i=2l'_{ij}$$
where $l'_{ij}$ is the $(i,j)$-entry of the matrix $L'=TLT^t$.
From the computations on the induced subalgebra we get that
$[A_{T}][A^{L}][A_{T^{-1}}]= [A^{L}(T)]=[A^{TLT^t}]$.

\begin{theorem}
The map
$$\Psi:(GL_n(k)/\Ent_2)\ltimes Sym_n(k) \rightarrow
BQ(k,E(n)), (T,L) \mapsto[A_{T^{-1}}][A^L]$$ is an injective
group homomorphism.
\end{theorem}
\pf It is not difficult to check that $\Psi$ is a group homomorphism. Let
$(\overline{T},L) \in (GL_n(k)/\Ent_2)\ltimes Sym_n(k)$ be such that
$\Psi[(\overline{T},L)]=[A^{L}][A_{{T}^{-1}}]=[k]$ in $BQ(k,E(n))$. Then
$[A_{T}]=[A^{L}]\in Sym_n(k)$. So, for every $N \in Sym_n(k)$,
$$[A^{TNT^t}]=[A_{T}][A^N][A_{T}]^{-1}=[A^{L}][A^N][A^{-L}]=[A^N].$$
Since the map $\chi$ in Proposition \ref{reduction} is injective
$TNT^t=N$ for all $N \in Sym_n(k)$, which implies that $T=\pm
Id$ and $\overline{T}=Id$. Then $[A_{T}]=[A^{L}]$ is trivial and $(L,\overline{T})=(0,\overline{Id}).$ \qed

\vspace{-0.5cm}
\section*{Acknowledgements}

The first author is partially supported by
Progetto Giovani Ricercatori number CPDG031245 of the University of
Padova. The second author is partially supported by project
BFM2002-02717 from Ministerio de Ciencia y Tecnolog\'{\i}a and
FEDER.

\end{document}